\documentclass[preprint,12pt]{elsarticle}
\usepackage{graphics}
\usepackage{subfigure}

\usepackage{latexsym}
\usepackage{amsmath}
\usepackage{graphicx}
\usepackage{times}
\usepackage{graphicx,multicol,enumerate,subfigure,amsmath}
\usepackage[usenames,dvipsnames,svgnames,table]{xcolor}
\usepackage{epstopdf}
\usepackage{fullpage}
\usepackage{setspace}
\usepackage[english]{babel}
\usepackage[version=3]{mhchem}
\usepackage{amsmath, amssymb, amsthm}
\usepackage{verbatim, latexsym}
\usepackage{colortbl}
\usepackage{float}

\usepackage{amsmath}
\usepackage{amsthm}
\usepackage{epsfig}
\usepackage{tikz}
\tikzset{>=latex}
\usepackage{setspace}
\usepackage{enumitem}
\usepackage{subfigure}
\usepackage{multirow}

\usepackage{soul}
\setstcolor{red}

\biboptions{sort&compress}
\journal{Journal of Computational and Applied Mathematics}

\onehalfspace

\begin{document}

\begin{frontmatter}

\title{ A mass 
conservative generalized multiscale finite element method applied to two-phase flow 
in heterogeneous porous media }
\author{\textbf{Michael Presho}$^{1*}$}
\cortext[cor1]{Email address : presho@ices.utexas.edu}

\author{\textbf{Juan Galvis}$^{2}$}

\address{$^{1}$ Institute for Computational and Engineering Sciences (ICES) \\
University of Texas at Austin \\
Austin, Texas}

\address{$^{2}$ Departamento de Matem\'{a}ticas \\
Universidad Nacional de Colombia \\
Bogot\'{a} D.C. }

\begin{abstract}
In this paper, we propose a method for the construction of locally conservative flux fields through a variation of the Generalized Multiscale Finite Element Method (GMsFEM). The flux values are obtained through the use of a Ritz formulation in which we augment
the resulting linear system of the continuous Galerkin (CG) formulation in the higher-order GMsFEM
approximation space. In particular, we impose the finite volume-based restrictions through incorporating a scalar Lagrange multiplier for each mass conservation constraint. This
approach can be equivalently viewed as a constraint minimization problem where
we minimize the energy functional of the equation restricted to the subspace of
functions that satisfy the desired conservation properties. To test the performance of the method we consider equations with heterogeneous permeability coefficients that have high-variation and discontinuities, and couple the resulting fluxes to a two-phase flow model. The increase in accuracy associated with the computation of the GMsFEM pressure solutions is inherited by the flux fields and saturation solutions, and is closely correlated to the size of the reduced-order systems. In particular, the addition of more basis functions to the enriched multiscale space produces solutions that more accurately capture the behavior of the fine scale model. A variety of numerical examples are offered to validate the performance of the method. 
\end{abstract}


\end{frontmatter}

\section{Introduction and problem statement}
In this paper, we primarily consider the equation given by

\begin{eqnarray}\label{eq:problem1}
-\mbox{div}(\Lambda k (x)\nabla p) = q \quad \mbox{in}~~\Omega \nonumber \\
p = p_D \quad \text{on}~~\Gamma_D  \\
-\Lambda k \nabla p \cdot \mathbf{n} = g_N \quad \text{on}~~\Gamma_N \nonumber 
\end{eqnarray}
where $k (x)$ is a heterogeneous field with high contrast. In particular,
we assume that there is a positive constant $k _{\min}$ such that $k (x)\geq k _{\min}>0$, while
$k (x)$ can have very large values (i.e., $k _{\max} / k _{\min}$ is large). Additionally, $\Lambda$ is a known mobility coefficient, $q$ denotes any external forcing, and $p$ is an unknown pressure field satisfying Dirichlet and Neumann boundary conditions given by $p_D$ and $g_N$, respectively. Here $\Omega$ a convex polygonal 
and two dimensional domain with boundary  $\partial \Omega= \Gamma_D\cup \Gamma_N$.

Let us consider a function in $H^1(D)$ whose trace on $\Gamma_D$ coincides with the 
given value $p_D$; we denote this function also by $p_D$. The
variational formulation of  problem  (\ref{eq:problem1}) is to find 
$p\in H^1(\Omega)$ with {$(p - p_D) \in H_D^1 = \{ w \in H^1(\Omega): ~ w|_{\Gamma_D} = 0 \}$} and such that 
\begin{equation}\label{eq:problem}
a(p,v)=F(v) - \langle g_N, v \rangle_{\Gamma_N} \quad \mbox{ for all } v\in H_D^1,
\end{equation}
where, for $p,v\in H^1(\Omega)$, the bilinear form $a$  is defined by
\begin{equation}\label{eq:def:a}
a(p,v)=\int_\Omega 
 \Lambda k (x)\nabla p(x)\nabla v(x) dx, 
\end{equation}
the functional $F$ is defined by 
\begin{equation}
{F(v)=\int_\Omega q(x)v(x)dx  }
\end{equation}
and the linear functional related to the boundary condition is given by
\begin{equation}
 \langle g_N, v \rangle_{\Gamma_N} =\int_{\Gamma_N} g_N(x) v(x) dl.
\end{equation}
A main goal of our work is to obtain conservative discretizations of the equations above. 
More specifically, construction of approximation that satisfy some given conservation of mass 
restriction on subdomains of interest. 
We note that a popular conservative discretization is the Finite Volume  (FV) method.  The classical 
FV  
discretization  provides and approximation 
of the solution in the space of piecewise linear functions  with respect to a triangulation while satisfying conservation of mass on elements of
 a dual triangualtion. 
When the approximation of the piecewise linear space is not 
enough for the problem at hand, advance approximation 
spaces need to be used. However, in some cases this requires a sacrifice of the conservation properties of the FV
method. In this work we present an extension of the FV method 
for general approximation spaces that enrich classical approximation 
spaces (such as the space of piecewise linear functions).
In particular, this conservative discretization can be used in conjunction with 
 recently introduced GMsFEM spaces. 

We note that FV methods that use higher degree piecewise polynomials have been introduced in the literature. The fact that the 
dimension of the approximation spaces is larger than the number of restrictions led the researchers of \cite{chen, chen2} to introduce additional 
control volumes to match the number of restrictions to the number of unknowns. 
An alternative approach is to consider a Petrov-Galerkin formulation with additional test functions rather that only piecewise constant functions on the dual grid. They where able to obtain stability of the method as well as error estimates. It is important to observe that the additional control volumes require additional 
computational work to be constructed and in some cases are not easy to construct (see also \cite{MR1906820, MR2113684}).
 Additionally, it is well know that 
piecewise smooth approximations spaces do not perform well for multiscale high-contrast problems
 \cite{egh12,jcp,eglp13oversampling, ge09_1,ge09_1reduceddim, eglw11, EGG_MultiscaleMOR, Review}. Another 
technique that one can use in order obtain a conservative method with richer approximations spaces is the following (see for instance \cite{CJimsfv} where 
the authors use a similar approach). In the discrete linear system obtained by a finite 
element discretization, it is possible to substitute appropriate number  of  equations 
by finite volume equations involving only the standard dual grid. This approach has the advantage that no additional 
control volume needs to be constructed. It may have the flexibility of both 
FV and FE procedures given as mass conservative fluxes and residual minimization properties.
Some previous numerical experiments suggest a drawback of this approach -   the resulting discrete problems may be 
ill-conditioned for large dimension coarse spaces, specially for higher order approximation spaces and multiscale problems. 

In this paper we propose the alternative of using a Ritz formulation and construct a solution procedure that 
combines a continuous Galerkin-type formulation that concurrently satisfies mass conservation restrictions. To this end,    
we augment the resulting linear system of the Galerkin formulation in the higher order approximation space to impose the finite volume restrictions. 
In particular,  
we do that by using a scalar Lagrange multiplier for each restriction. This approach can be equivalently viewed as a constraint minimization 
problem where we minimize the energy functional of the equation restricted to the subspace of functions that satisfy the 
conservation of mass restrictions. Then, in the Ritz sense, the obtained solution
is the best among all functions that satisfy the mass conservation restriction.

As a main application of the 
techniques presented here, we  consider the case where the coefficient $k $
has high-variation and discontinuities (not necessarily aligned with the coarse grid). 
For this problem it is known that higher order approximation is needed. Indeed, 
in some cases robust approximation properties, independent of the contrast, are required. See for instance 
\cite{ ge09_1,ge09_1reduceddim, eglw11} where it is demonstrated that classical multiscale methods (\cite{eh09}) 
do not render robust approximation properties in terms of the contrast. It is shown that 
one basis functions per coarse  node (with the usual support) is not enough 
to construct adequate coarse spaces \cite{ge09_1reduceddim, MR2861243}. A similar issue can be expected for 
the multiscale finite volume method developed in \cite{MR2575046,MR3109801,MR2123109} and related works, when applied to 
high-contras multiscale problems since the approximation spaces have similar 
approximation properties. In the case of Galerkin formulations, 
robust approximation properties are obtained by using the 
Generalized Multiscale Finite Element  Method (GMsFEM) framework.
The main goal of GMsFEMs is to construct  coarse spaces for Multiscale Finite Element Methods (MsFEMs)
that result in accurate coarse-scale solutions. 
This methodology was first developed in \cite{egh12,jcp,eglp13oversampling} 
based on some previous works \cite{ge09_1,ge09_1reduceddim, eglw11, EGG_MultiscaleMOR, Review}. 
A main ingredient in the construction is the use of an approximation of  local eigenvectors
 (of carefully selected local eigenvalues problem)  to construct the coarse spaces.
Instead of  using one coarse function per coarse node as in classical MsFEM, in the GMsFEM
it was proposed to use several multiscale basis functions per coarse node.
These basis functions represent
important features of the solution within a coarse-grid 
block and they are 
computed using eigenvectors of an eigenvalue problem.
For applications to high-contras problems, methodologies to keep small the dimension of the resulting coarse space where successfully proposed 
(\cite{ge09_1reduceddim}).
That paper made use of coarse spaces that somehow incorporated important modes of a (local) energy 
related to the problem motivated the general version of the GMsFEM.
The much more general GMsFEM 
was then developed in \cite{egh12} where several more practical options to 
compute important modes to be include in the coarse space were used; see also \cite{egt11} for an earlier construction.
It is important 
to mention that the methodology 
in \cite{egh12} was designed for parametric and nonlinear problems and 
can be applied in variety of applications, although a more extensive review of such developments is not contained herein.

An important consideration of the proposed method is that the GMsFEM methodology does not guarantee conservation of mass properties 
such as those of Finite Volume formulations. In this paper 
we design a mass conservative GMsFEM method. In particular, we impose the 
mass conservation constraints over control volumes by using Lagrange multipliers.
In doing so, we numerically validate that the convergence properties of the GMsFEM are 
maintained while the approximate solution simultaneously satisfies the required conservation properties. We mention that in \cite{MR3133358} 
some successful applications of  GMsFEM to two-phase problems were
presented. In that work, the authors developed a 
technique in which GMsFEM solutions could be post-processed to yield mass conservative fluxes for use in two-phase modeling. While the motivation of the previous work in
 \cite{MR3133358} is similar, the current technique yields a system that requires no post-processing, and a solution that automatically yields conservation. In particular, in the present
  work no additional equations must be solved (in notable constrast to \cite{MR3133358}) in order to ensure coarse-grid conservation, as a single global solve automatically yields the
   desired properties. As a result, no additional computational resources must be allocated, and issues such as ill-conditioning of the localized systems may be circumvented. Additionally, the proposed method ensures that the solution obtained through the method proposed in this work is the best among all functions (in the Ritz sense) satisfying the mass conservation restrictions.

Finally we mention the works on MsFV. When dealing with multiscale problems, the lack of accuracy in the use of 
classical spaces had led researchers to introduce enriched spaces (as well as 
iterative procedures).  In \cite{CJimsfv} the authors construct an enrichment of the 
initial coarse space spanned by nodal basis functions. We also note that a hybrid 
finite-volume/Galerkin formulation for the coarse-scale problem is devised. The authors show that the resulting method has all 
features of the MsFV method, but is more robust and shows improved convergence 
properties if used in an iterative procedure. In comparison, in this work we use a 
GMsFEM framework were the additional degrees of freedom are related to local eigenvalue problems. This construction is related to the approximation properties 
the resulting space will have and is motivated by the numerical analysis of the 
interpolation operator into the coarse space. Furthermore, instead of 
a hybrid finite-volume/Galerkin formulation, we employ a Ritz formulation in the space of functions satisfying the conservation of mass.

The rest of the paper is organized as follows. 
In Section \ref{lagrange}, we recall how to impose restrictions 
using Lagrange multipliers in general and, in particular, for the 
mass conservation restrictions. In Section \ref{discrete} we
formulate a conservative coarse problem that follows the framework of Section \ref{lagrange}. Section \ref{smooth} 
is dedicated to present some relevant discussions for the case of smooth 
coefficients and solutions. In Section \ref{tpmodel} we present the 
two-phase flow model problem and the overall solution algorithm. 
In Section \ref{gmsfem} we summarize some basic construction topics associated with the 
GMsFEM methodology. In Section \ref{numerical} we present some representative numerical results to illustrate the successful performance of our method for the 
case of single- and two-phase flow problems in high-contrast heterogeneous porous media. Some concluding remarks are finally offered in Section \ref{conclusion}.

\section{Linear restrictions using Lagrange Multipliers}
\label{lagrange}

Problem (\ref{eq:problem}) is equivalent to the minimization 
problem: Find $p$ with  $(p - p_D) \in H_D^1$ and such that 
\begin{equation}
p  =\arg\min_{v} \mathcal{J}(v)
\end{equation}
where the minimum is taken over $v$ such that $v-p_D\in H^1_D$ and 
\begin{equation}\label{eq:def:J}
\mathcal{J}(v)=\frac{1}{2}a(v,v)-F(v)+\langle g_N, v \rangle_{\Gamma_N}.\end{equation} 

Let $p$ be the solution of \eqref{eq:problem} and   $\tau_i$, $i=1,\dots,M$, be $M$ continuous linear functionals 
on $H^1(D)$. Define $m_i=\tau_i(p)$, $i=1,2,\dots,M$. The problem above is 
equivalent to: Find $(p-p_D)\in H_D^1$ such that 
\begin{equation}\label{eq:problem-with-restriction}
p =\arg\min_{v\in \mathcal{W}} \mathcal{J}(v)
\end{equation}
where 
\[
\mathcal{W}=\{ v: v-p_D\in  H^1_D \mbox{ and } \tau_i(v)=m_i, \quad i=1,\dots,M\}.
\]
Problem (\ref{eq:problem-with-restriction}) above 
can be view as Lagrange multipliers min-max optimization 
problem. See \cite{MR2168342} and references therein. 

Then, in case an approximation of $p$, say $p^h$, it is required to satisfy the constraints  $\tau_i(p^h)=m_i$, $i=1,2,\dots,M$, we can 
discretize directly the formulation \eqref{eq:problem-with-restriction}. In particular, we can 
apply this approach to a set of mass conservation
restrictions used in finite volume discretizations.

\subsection{Mass conservation in fine control volumes} \label{masscv}
In order to ensure fine-scale conservation, we start by selecting  control volumes   
$\{ V_{i,f} \}_{i=1}^{M_f}.$ We assume that each $V_{i,f}$ is a subdomain of $\Omega$ with polygonal
boundary and $i=1, \dots,M_f$. 
If $q\in L^2$ we have that (\ref{eq:problem}) is equivalent
to:
Find  $p$ with {$(p-p_D)\in H_D^1$} and such that 
\begin{equation}\label{eq:problem-with-restriction2}
p=\arg\min_{v\in \mathcal{W}} \mathcal{J}(v)
\end{equation}
where the subset of functions that satisfy the mass conservation restrictions is defined by
\[
\mathcal{W}=\left\{ v\in H^1_D: \int_{\partial V_{i,f}}
- \Lambda k  \nabla v\cdot\mathbf{n} =\int_{V_{i,f}}q \quad \text{ for all } V_{i,f} \right\}.
\]

The Lagrange multiplier formulation of problem (\ref{eq:problem-with-restriction}) 
can be written as:
Find $p$ with $(p - p_D) \in H^1_D$ and $\lambda \in\mathbb{R}^{M_f}$ that solves,
\begin{equation}\label{eq:discrete-problem-with-restriction-lag}
\max_{\mu\in \mathbb{R}^{M_f}}\min_{v\in H^1_D} \mathcal{J}(v)-
(\overline{a}(p,\mu)-\overline{F}(\mu)).
\end{equation}
Here, the average flux  bilinear form $\overline{a}: H^1_D\times \mathbb{R}^{M_f}\to \mathbb{R}$ is defined by
\begin{equation}\label{eq:def:overline-a}
\overline{a}(v,\mu)=\sum_{i=1}^{M_f}
\mu_i\int_{\partial V_{i,f}} - \Lambda k \nabla v\cdot \mathbf{n} \quad \mbox{ for all }  v\in H_D^1 \mbox{ and } \mu\in \mathbb{R}^{M_f}.
\end{equation}
The functional $\overline{F}:\mathbb{R}^{M_f} \to \mathbb{R}$ is defined by\[\overline{F}(\mu)=\sum_{i=1}^{M_f} \mu_i\int_{V_{i,f}}q\quad  \mbox{ for all } \mu \in \mathbb{R}^{M_f}.\]
The first order conditions of the min-max problem above  give the following saddle point problem:
Find $p$ with $(p-p_D)\in H^1_D$ and $\lambda\in\mathbb{R}^{M_f}$ that solves,
\begin{equation}\label{eq:saddlepoint}
\begin{array}{llr}
a(p,v)+\overline{a}(v,\lambda)&=F(v)-\langle g_N, v \rangle_{\Gamma_N} &\mbox{ for all } v\in H^1_D, \\
\overline{a}(p,\mu) &=\overline{F}(\mu)& \mbox{ for all } \mu \in \mathbb{R}^{M_f}.\\
\end{array}
\end{equation}
See for instance \cite{MR2168342}.

\section{Conservative discrete coarse problem}
\label{discrete}

Let $ H$ be a coarse-mesh parameter and 
$\mathcal{T}^H$ be a coarse-scale triangulation. We assume that $H$ does not necessary 
resolve all  the variation of the coefficient.
In what follows we introduce the finite element space $V^H$ associated 
with the coarse resolution $H$.  In applications concerning heterogeneous multiscale media, the standard finite element spaces on $\mathcal{T}^H$  do not offer good approximation properties  and therefore some special enrichment of this space is needed in order to achieve some acceptable approximation properties.

\begin{figure}[tb]
 \centering
   \includegraphics[width = 0.95\textwidth, keepaspectratio = true]{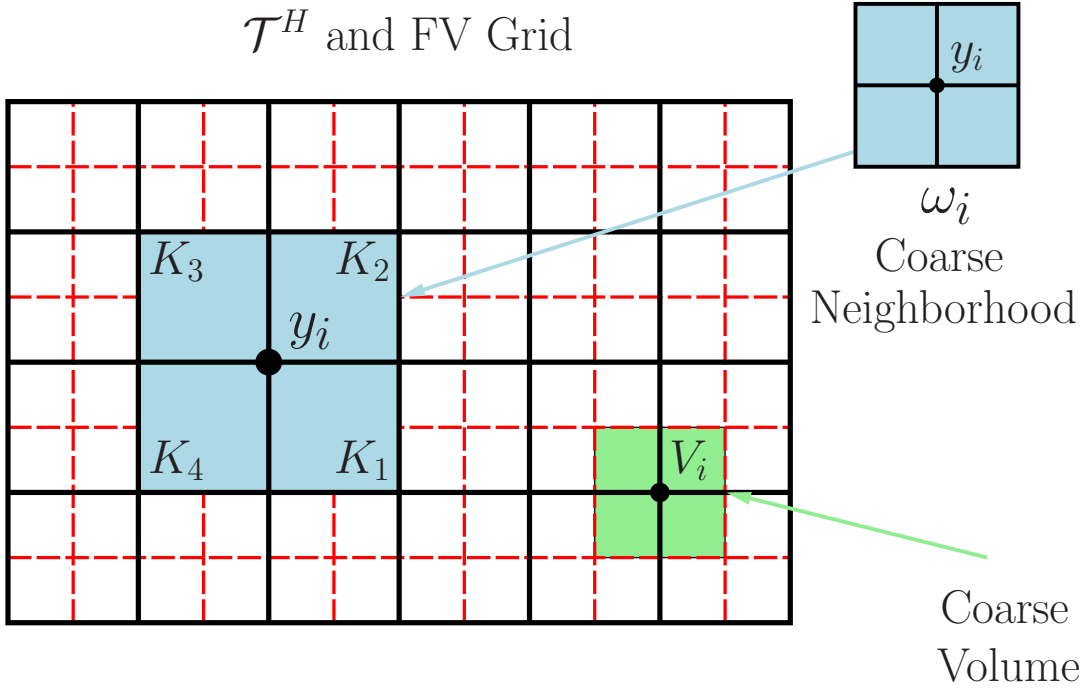}
 \caption{ Illustration of a coarse discretization depicting coarse neighborhoods and control volumes }
 \label{schematic}
\end{figure}

We denote by $\{y_i\}_{i=1}^{N_v}$ the vertices of the coarse mesh
$\mathcal{T}^H$ and define the neighborhood of the 
node $y_i$ by
\begin{equation}\label{eq:def:omegai}
\omega_i=\bigcup\{ K_j\in\mathcal{T}^H; ~~~ y_i\in \overline{K}_j\}.
\end{equation}
%
See Fig. \ref{schematic} for an illustration of the coarse-scale discretization depicting coarse neighborhoods and control volumes. 

Using the coarse mesh $\mathcal{T}^H$ 
we introduce  a set of coarse
 basis functions $\{\Phi_i\}$. The basis functions
are supported in $\omega_i$; however, for one specific neighborhood $\omega_i$, there
may be multiple basis functions.
We define the associated coarse space by 
\begin{equation}\label{eq:def:V0}
V_0=\mbox{span} \{ \Phi_i\}_{i=1}^{\text{dim}(V_0)}.
\end{equation}
Let $p_D^H$ be a discrete interpolation of the boundary data $p_D$. 
The classical multiscale solution of (\ref{eq:problem})   is given by $p_{ms}$ 
with $p_{ms}-p_D^H\in V_0$ such that 
\begin{equation}\label{eq:ritzgalerkin}
p_{ms}=\arg\min \mathcal{J}(v)
\end{equation}
where the minimum is taken over the set of $v$ with $v-p^H_D\in V_0$.
%
The equivalent  matrix form of the Ritz-Galerkin formulation (\ref{eq:ritzgalerkin}) is given by
\begin{equation}\label{eq:matrixproblem}
A_0p_{ms}=b
\end{equation}
where the coarse matrix $A_0$ and the coarse load vector $b$ are such that for all $u,v\in V_0$ we have 
\begin{equation}\label{eq:def:A0}
u^TA_0v=a(u,v) 
\quad \mbox{ and }\quad   v^Tb=F(v) - \langle g_N, v \rangle_{\Gamma_N} .
\end{equation}

The multiscale solution $p_{ms}$ of problem (\ref{eq:matrixproblem}) does not necessarily have any conservation
of mass property. Therefore, the classical multiscale 
solution might not be appropriate form some applications where it is imperative that 
approximations have some conservation properties in the form of average flux 
through some control volumes.
To obtain multiscale solutions with the required conservation properties one can proceed as described in 
Section \ref{lagrange}.
To that end, we let $V_0$ be defined as before and consider the set of all discrete functions that 
satisfy the required approximation properties,
\begin{equation*} \label{coarseconserv}
\mathcal{W}_{0}=\{ v: v-p_D^H\in V_0 \mbox{ and } \int_{\partial V_{i,c}}
-\Lambda k  \nabla v\cdot\eta =\int_{V_{i,c}}q,  \quad i=1,\dots,M_c\},
\end{equation*}
where $V_{i,c}$ denotes a coarse dual-grid volume. 

We then consider the following discrete formulation that takes into account the required restrictions. The approximation of
the solution of (\ref{eq:problem1}) is to find $p_{fv}\in \mathcal{W}_0$  such that 
\begin{equation}\label{eq:discrete-problem-with-restriction}
p_{fv}=\arg\min_{v\in \mathcal{W}_0} \mathcal{J}(v).
\end{equation}
This is a minimization problem with linear constraints 
similar to  (\ref{eq:problem-with-restriction}) - Indeed, it is the associated Ritz-Galerkin formulation; see Section \ref{lagrange}. 
This problem can be equivalently written in the form analogous to (\ref{eq:discrete-problem-with-restriction-lag}) or (\ref{eq:saddlepoint}).
The matrix  formulation stating 
the first order conditions of this minimization problem is written as
\begin{equation}\label{eq:matrix-fv}
\left[\begin{array}{cc}
A_0&\overline{A}^T\\
\overline{A}& 0\\
\end{array}\right]\left[\begin{array}{c}p_{fv}\\\lambda\end{array}\right]=
\left[\begin{array}{c}b\\\overline{b}\end{array}\right],
\end{equation}
where the matrix $A_0$ and the vector $b$ are defined in 
(\ref{eq:def:A0}) and the finite volume matrix $\overline{A}$ is defined by
\[
\mu^T \overline{A} v=\overline{a}(v,\mu)\quad \mbox{ and }
 \quad \mu^T\overline{b}=\overline{F}(\mu) \quad \mbox{ for all } 
 \mu \in \mathbb{R}^{M_c} \mbox{ and } v\in H_0^1(D).
\]

\section{The case of smooth solutions}\label{smooth}

In this section, in order to motivate the use of the the formulation \eqref{eq:discrete-problem-with-restriction} for multiscale problems, we 
make some comments on the method described so far for  the case of smooth coefficients and  higher-order degree polynomials for finite volume methods.

As mentioned in the introduction, when higher order degree polynomials are used within  a finite volume framework 
one can proceed in different ways. Let us assume, only for this section and with motivation purposes, 
that  $V_0$ is the space of piecewise polynomials  $P^r(\mathcal{T}^{H})$ with $r\geq 1$. Let us then introduce
 the dual grid $\mathcal{T}^{H*}$. The dimension of the approximation space $V_0$ is higher that the number of 
controls volumes or dual elements. The number of control volumes matches the dimension of the subspace 
of piecewise linear functions $P^1(\mathcal{T}^{H})$  contained in $V_0$. When $r>1$ and if we use only piecewise
constant functions on the dual grid as a test function we do not obtain a square system.  In order to obtain a square 
linear system one have to introduce additional test functions. One can proceed as follows.  
\begin{enumerate}
\item Construct additional control volumes and test the approximation spaces against piecewise constant functions 
over the total of control volumes (that include the dual grid element plus the additional control volumes). 
We mention that constructing additional control volumes is not an easy task and might be computationally expensive.
 We refer the interested reader to \cite{chen, chen2015,  chen2} for additional details.
\item Use as additional the basis functions the basis functions that correspond to nodes other than vertices to obtain a FV/Galerkin formulation.
This option has the advantage that no geometrical constructions  have to be carried out. On the other hand, 
this formulation seems difficult to analyze. Also, some preliminary numerical tests suggest that the resulting linear
 system becomes unstable for higher order approximation spaces (especially for the case of high-contrast multiscale coefficients).
\item Use the Ritz formulation with restrictions \eqref{eq:discrete-problem-with-restriction}.
\end{enumerate}

Note that if $r=1$, in the
linear system associated to \eqref{eq:matrix-fv}, the restriction matrix corresponds to the usual 
finite volume matrix.  This matrix is known to be invertible. In this case, the affine space 
$\mathcal{W}_0$ is a singleton. 
Moreover, the only  function $u$ satisfying the restriction is given by $u=(\overline{A})^{-1} \overline{b}$. 
The Ritz formulation \eqref{eq:discrete-problem-with-restriction} reduces to the classical finite volume 
method.

If we switch to $r>1$, the dimension of the space $\mathcal{W}_0$ can only increase.  Therefore and due to ellipticity of the energy functional, the Ritz formulation with restrictions 
\eqref{eq:discrete-problem-with-restriction} makes sense and it will give the best approximation of the 
solution in the space of functions that satisfy the restrictions. Then, in the Ritz sense, the solution of 
\eqref{eq:discrete-problem-with-restriction} is not worse that any of the solutions obtained by 
the methods 1. or 2. mentioned above. Furthermore, the solution of the associated linear system 
\eqref{eq:matrix-fv}, which is a saddle point linear system, can be readily implemented using
efficient solvers for the matrix $A$ (or efficient solvers for the classical finite volume matrix $\overline{A}$);
See for instance \cite{MR2168342}.
Additionally, we mention that the analysis of the method can be carried out using usual tools for the analysis of 
restricted minimization of energy functionals and mixed finite element methods. The numerical analysis of our methodology 
is under current investigation and it will presented elsewhere.\\

The comments and observations of this section are the main motivation for the methodology developed here to target 
fluid flow in porous media where, as it is well known, advanced and sophisticated approximation spaces have to be used. 
Furthermore, in some applications such as multiphase flow, conservation of mass is a main requirement.

\section{Two-phase model problem} \label{tpmodel}
We emphasize that solving the pressure equation in \eqref{eq:problem1} is done within the context of a two-phase flow model. 
In particular, we are interested in treating a problem confined to a domain 
$\Omega$ in which the subsurface is assumed to exhibit high-contrast features.
 The heterogeneous reservoir is equipped with a well in which water is injected to displace the trapped oil towards the production wells. 
 The water and oil phases (which we denote by $o$ and $w$, respectively) are assumed to be immiscible, and we consider a gravity-free
  enviroment in which the pore space of the reservoir is fully saturated. Additionally, we assume that any capillary effects are negligible.
   Under such assumptions, combining Darcy's law 
  with a statement of conservation of mass yields governing equations of the form 

\begin{equation} \label{pres}
\nabla \cdot \mathbf{v} = q,
\end{equation}
where the total Darcy velocity is given by

\begin{equation} \label{vel}
\mathbf{v} = -\Lambda(S) k(x) \nabla p.
\end{equation}
The pressure equations given by \eqref{pres} is coupled to a transport equation 

\begin{equation} \label{sat}
\frac{\partial S}{\partial t} + \nabla \cdot \left( \mathbf{v} f(S)  \right) = q_w,
\end{equation}
where $S$ denotes the water saturation, and $q, q_w$ denote any external forcing. The total mobility coefficient $\Lambda(S)$ and flux function $f(S)$ are respectively given by 

\begin{equation*}
\Lambda(S) = \frac{k_{rw}(S)}{\mu_w} + \frac{k_{ro}(S)}{\mu_o} \quad \text{and} \quad 
f(S) = \frac{k_{rw}(S) / \mu_w}{\lambda(S)},
\end{equation*}
where $k_{r \alpha}, \alpha = w,o$ is the relative permeability of the $\alpha$ fluid phase. 

\subsection{Solution algorithm} \label{solalg}
\begin{figure}[tb]
 \centering
   \includegraphics[width = 1.0\textwidth, keepaspectratio = true]{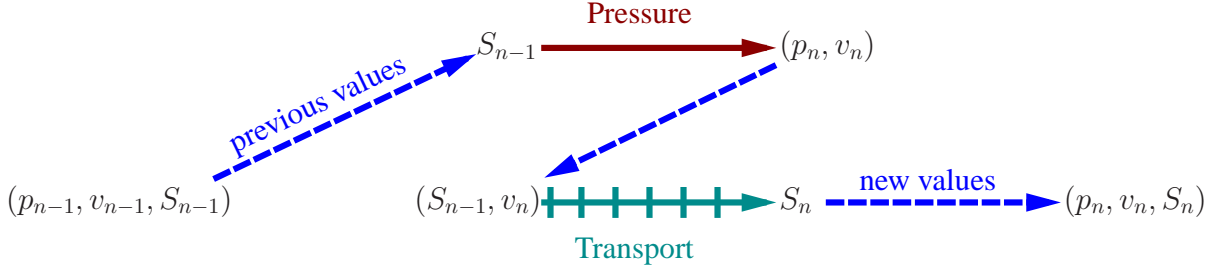}
 \caption{ An illustration of the operator splitting technique used to solve the two-phase model problem }
 \label{opsplit}
\end{figure}

In order to solve the coupled two-phase model given in Eqs. \eqref{pres}, \eqref{vel}, and \eqref{sat} we employ an operator splitting 
technique where the saturation from the previous time step is used to update the mobility coefficient required of the pressure solve and 
subsequent velocity calculation (see, e.g., \cite{as79}). Once this velocity is available, it is used in conjunction with an explicit saturation 
time-marching scheme for a specified number of time steps. The updated saturation is used again in order to update the pressure,
 and the process is continued until a final simulation time is reached. See Fig. \ref{opsplit} for a schematic of the the operator splitting. 

To discretize the saturation equation, we integrate Eq. \eqref{sat} with respect to time, and then over some volume $V_i \in \Omega$. Applying the left end-point quadrature rule to a second term in time, and integrating by parts yields the following expression

\begin{equation*}
\text{meas}(V_i) (S_{i,n} - S_{i, n-1}) + \Delta t \int_{\partial V_i} \mathbf{v} \cdot \mathbf{n} f(S_{i,n-1})  = \int_{V_i} q_w ,
\end{equation*}
where the errors terms have been neglected and

$$
S_{i,n} \approx \frac{1}{\text{meas}(V_i) } \int_{V_i} S(x, t_n).
$$
In addition to the explicit time stepping in the above expression, we apply a non-oscillatory upwinding scheme on the flux term $\int_{\partial V_i} \mathbf{v} \cdot \mathbf{n} f(S_{i,n-1}) $ (see, e.g., \cite{t99} for a review of upwinding schemes on rectangular meshes). As mentioned in Section \ref{lagrange}, it is crucial that that the numerical flux approximation $\mathbf{v}_h$ satisfies a local conservation property. In particular, we wish to seek discrete solutions that satisfy 

\begin{equation} \label{conserv}
\int_{\partial V_i} \mathbf{v}_h \cdot \mathbf{n} = \int_{V_i} q \quad \text{for each} ~~ V_i.
\end{equation}
We emphasize that a main consideration of this paper is to ensure the conservation property in \eqref{conserv} through the use of Lagrange multiplier restrictions (refer to Section \ref{lagrange}) within a suitable coarse-scale solution space (refer to Section \ref{discrete}). In turn, in the next section we describe the systematic construction of the Generalized Multiscale Finite Element Method (GMsFEM) coarse solution space.

\section{Generalized multiscale finite element method}
\label{gmsfem}
In this section we focus on high-contrast multiscale problems and summarize a GMsFEM construction 
of $V_0$ used in Section \ref{discrete}. For a more detailed description of the development of the GMsFEM methodology, see \cite{egh12,jcp,eglp13oversampling} and references therein.

\subsection{Spectral enrichment} 
We start by chosing and initial set of basis functions that form a partition of unity. 
The space generated by this basis functions is enriched using a local spectral problem.
We use the multiscale basis functions partition of unity  with linear
boundary conditions (\cite{eh09}). We have one function per coarse-node and it is defined by 
\begin{eqnarray}\label{eq:MsFEM_standard1}
-\mbox{div}(k  \nabla {\chi}_{i} )&=&0\quad \mbox{for}\ K \in \omega_i \\
 {\chi}_{i}&=& {\chi}_{i}^0  \quad 
\mbox{ on } \partial K, \nonumber
\end{eqnarray}
where ${\chi}_{i}^{0}$ is a standard linear partition of unity function.

For each coarse node neighborhood $\omega_i$, consider the eigenvalue problem
\begin{equation}
\label{eq:eig:prob}
-\mbox{div}(k  \nabla \psi_\ell^{\omega_i})=\sigma_\ell^{\omega_i} \widetilde{k } \psi_\ell^{\omega_i},
\end{equation}
with homogeneous Neumann
boundary condition on $\partial \omega_i$. Here $\sigma_\ell^{\omega_i}$
 and  $\psi_\ell^{\omega_i}$ are eigenvalues and
eigenvectors in $\omega_i$ and
$\widetilde{k }$ is defined by
\begin{equation}\label{eq:def:kappatilde}
\widetilde{k }=
k \sum_{j=1}^{N_v}
{H^2}|\nabla \chi_j|^2. 
\end{equation}
We use an ascending ordering on the eigenvectors, 
$
\sigma_1^{\omega_i}\leq \sigma_2^{\omega_i}\leq....
$\\

Using the partition of unity functions from Eq. \eqref{eq:MsFEM_standard1} and eigenfunctions from Eq. \eqref{eq:eig:prob}, we then construct a set of enriched multiscale basis functions given by
$\chi_i \psi^{\omega_i}_\ell$ for selected eigenvectors $\psi^{\omega_i}_\ell$. Using
$L_i$ to denote the number of basis functions from the coarse region $\omega_i$, we then define the coarse GMsFEM space by
\[V_0=\mbox{span}\{\Phi_{i,\ell}=\chi_i \psi_\ell^{\omega_i},\quad i=1,\dots,N_v,\quad \ell=1,\dots,L_i\}. \]
For more details, motivation of the construction, and approximation properties of the space $V_0$ as well as the choice of the initial partition of unity basis functions we refer the interested reader to \cite{egh12,jcp,eglp13oversampling}.

Summarizing, in order to solve problem (\ref{eq:problem1}) for the pressure we use the GMsFEM coarse space $V_0$ constructed in 
this section. Additionally, in order to obtain conservative solutions with respect to a dual coarse grid, we use the discretization presented in Section \ref{discrete}. In particular, we solve the saddle point problem (\ref{eq:matrix-fv}) using the appropriate approximation spaces and conservation constraints. We also recall Section \ref{tpmodel} for a description the overall 
solution algorithm used to solve the two-phase flow problem.

\subsection{A flux downscaling procedure} \label{downscale}
At this point we carefully distinguish the scales at which we wish to calculate respective flux values. Recalling the discussions from Section \ref{masscv} and Section \ref{discrete}, we 
note that sets of control volumes $\{ V_{i,f} \}_{i=1}^{M_f}$ or $\{ V_{i,c} \}_{i=1}^{M_c}$ are selected in order to incorporate the associated mass conservation restrictions of the 
fine and coarse problems, repsectively. In the case of the fine-scale problem, we simply assume that the set of control volumes $\{ V_{i,f} \}$ coincides with the nodal values of the 
mesh. In particular, the conservation constraints are taken directly on the fine scale. The case when we use GMsFEM is similar in the sense that we initially impose the constraints on the 
(larger) $\{ V_{i,c} \}$  volumes associated with the coarse mesh discretization (see Fig. \ref{schematic}). As such, a global GMsFEM-FV solve directly yields flux parameters that satisfy 
a coarse analogue of conservation. While this level of conservation may hold for some target applications (see, e.g., \cite{gppw11}), in this paper we wish to construct fluxes that are 
conservative directly on the fine scale for a direct means of comparison with the fine-scale FV technique. 

In constructing fine-scale mass conservative fluxes, we recall that the coarse solve from Section \ref{discrete} yields the coarse conservation property

$$
\int_{\partial V_{i,c}} \mathbf{v}_H \cdot \mathbf{n} = \int_{V_{i,c} } q  \quad \text{for all}~~ V_{i,c},
$$
which can analogously be thought of as the compatibility condition in $V_{i,c}$. As such, we may formulate a fully Neumann boundary value problem 

\begin{eqnarray} \label{dscale}
&-\nabla \cdot \left( \Lambda k \nabla p_{V_{i,c}}     \right) = q \quad \text{in}~~ V_{i,c} \\
&-\Lambda k \nabla p_{V_{i,c}} \cdot \mathbf{n}  = \mathbf{v}_H \cdot \mathbf{n} \quad \text{on}~~ \partial V_{i,c} \nonumber, 
\end{eqnarray}
where the known $\mathbf{v}_H = -\Lambda k \nabla p_{fv}$ is evaluated pointwise on the boundaries $\partial V_{i,c}$. In particular, after obtaining the coarse solution $p_{fv}$ 
and the coarse-scale conservative flux $\mathbf{v}_H$, we solve the set of localized problems in Eq. \eqref{dscale} for every $V_{i,c}$ in $\Omega$ using any method that produces 
the desired fine-scale conservation. For consistency within this paper, we use the fine-scale FV technique to ensure fine conservation. A hallmark advantage of the downscaled post-
processing procedure from \eqref{dscale} is that the localized problems are independent from one another. In particular, the problems are naturally parallelizable (as should also be 
carefully noted for Eqs. \eqref{eq:MsFEM_standard1} and \eqref{eq:eig:prob}), and may be independently distributed to CPU and/or GPU multi-core clusters for a significant gain in 
computational efficiency \cite{hw97}.

 \section{Numerical results}
\label{numerical}

In this section we offer a variety of numerical examples to test the performance of the method introduced in Section \ref{lagrange}. In particular, we solve the model problem in Eq. 
\eqref{eq:problem1} using the constrained discretization techniques from Sections \ref{discrete} and \ref{gmsfem}. A main goal is to address the accuracy associated with the reduced-
order  conservative GMsFEM approach as compared to the fine-scale FV approach. 

\subsection{Single-phase pressure} \label{spresult}
\begin{figure}[tb]
 \centering
   \includegraphics[width = 0.5\textwidth, keepaspectratio = true]{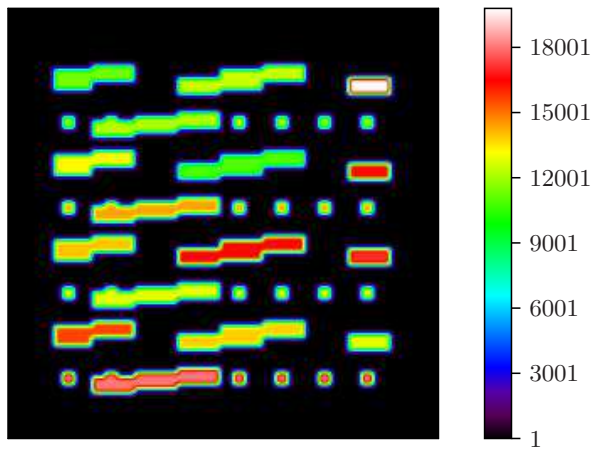}
 \caption{ A high-contrast permeability coefficient with inclusions and channels.  }
 \label{coeff}
\end{figure}

For the first set of examples, we employ both the fine-scale FV and GMsFEM-FV approaches to solve Eq. \eqref{eq:problem1}. Throughout the section we consider solutions that are 
obtained through solving the equation on the unit square domain $D = [0,1] \times [0,1]$. We impose boundary conditions of $p_L = 1$ and $p_R = 0$ on the left and right 
boundaries of the domain, along with no-flow (i.e., zero Neumann) conditions on the top and bottom. The fine coefficient (and reference solution) are posed on a $100 \times 100$ fine 
mesh that yields a global system of size $N_f = 20402$. See Fig. \ref{coeff} for an illustration of high-constrast structure that is considered in this section. More specifically, the system 
from Eq. \eqref{eq:matrix-fv} is of size $N_f = \text{dim}(V^h) + M_{f} = 10201 + 10201$, where $\text{dim}(V^h)$ denotes the dimension of the space in which the fine-grid 
pressure solution is represented, and $M_{f}$ denotes the number of dual-grid volumes where the finite volume constraints are imposed. In constructing the coarse-grid solutions in 
this section, we consider a $10 \times 10$ coarse mesh with varying levels of constrained spectral enrichment. In particular, we obtain systems of size $N_c = \text{dim}(V_0) + M_c$ 
where $\text{dim}(V_0)$ denotes the number of degrees of freedom associated with the coarse enrichment, and $M_c$ denotes the number of coarse dual-grid volumes where the 
associated finite volume constraints are imposed. For this set of examples we emphasize that the number of coarse dual-grid volume constraints is fixed at $M_c = 121$, which 
corresponds to the number of coarse nodal values (and surrounding volumes) of the domain (cf. Fig. \ref{schematic}). 

\begin{figure}[tb]
 \centering
   \hspace*{-0.04\textwidth}
  \subfigure[Fine $N_f = 20402$]{\label{fig:solf}
    \includegraphics[width = 0.33\textwidth, keepaspectratio = true]{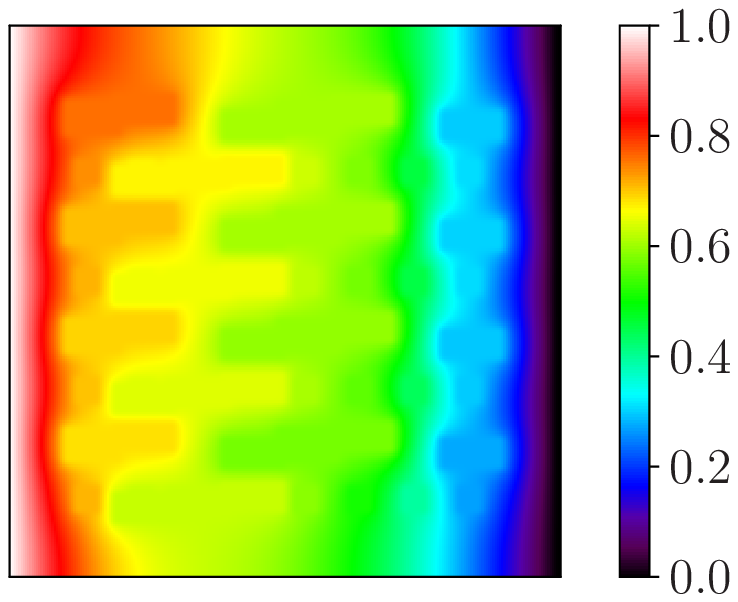}
    }
    \hspace*{-0.04\textwidth}
  \subfigure[Coarse $N_c = 323$]{\label{fig:solc2}
     \includegraphics[width = 0.33\textwidth, keepaspectratio = true]{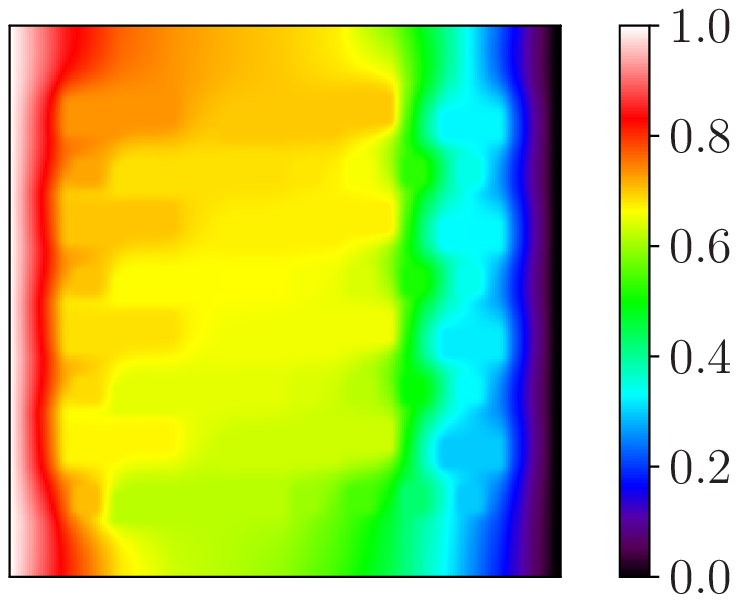}
    }
    \hspace*{-0.04\textwidth}
  \subfigure[Coarse $N_c = 647$]{\label{fig:solc6}
     \includegraphics[width = 0.33\textwidth, keepaspectratio = true]{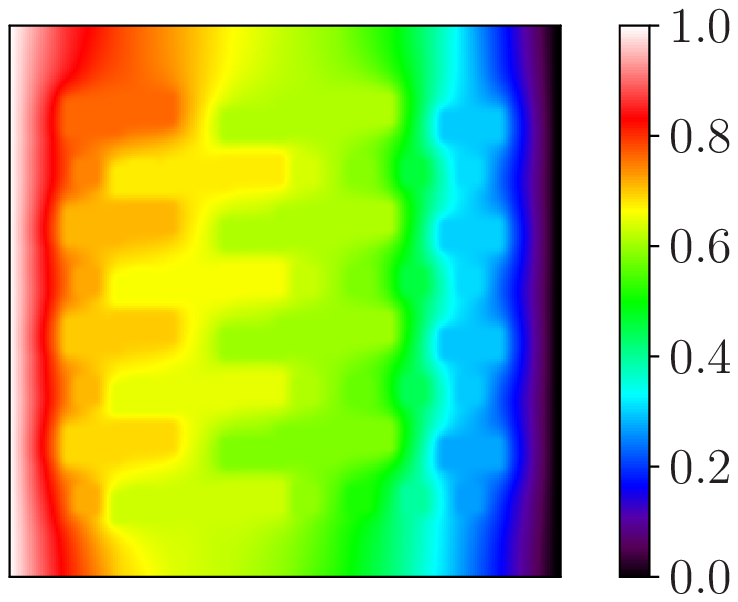}
    }
  
 \caption{ A comparison between the fine pressure solution and increasingly accurate coarse-grid  pressure solutions. }
 \label{solplot}
\end{figure}

As an initial motivation, we offer a comparison of single-phase pressure solutions in
 Fig. \ref{solplot}. The coarse solutions in Figs. \ref{fig:solc2} and \ref{fig:solc6} 
 were obtained through solving the global equation using two levels of coarse mesh 
 enrichment. We can see from Fig. \ref{fig:solc2} that a system of size $N_c = 323$ 
 offers a somewhat crude approximation to the reference solution, yet that a system 
 of size $N_c = 647$ (cf. Fig. \ref{fig:solc6}) yields a solution that is nearly 
 indistinguishable from the reference solution. In either case, we emphasize that the 
 coarse systems are much smaller than the system of size $N_f = 20402$ that is used 
 to obtain the fine-grid reference solution.  

\begin{figure}[tb]
 \centering
   \hspace*{-0.04\textwidth}
  \subfigure[Fine $N_f = 20402$]{\label{fig:fluxf}
    \includegraphics[width = 0.33\textwidth, keepaspectratio = true]{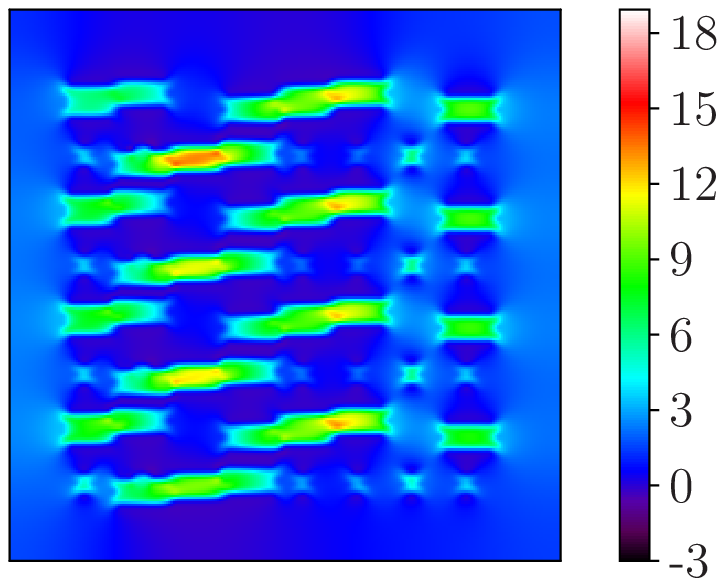}
    }
    \hspace*{-0.04\textwidth}
  \subfigure[Coarse $N_c = 323$]{\label{fig:fluxc2}
     \includegraphics[width = 0.33\textwidth, keepaspectratio = true]{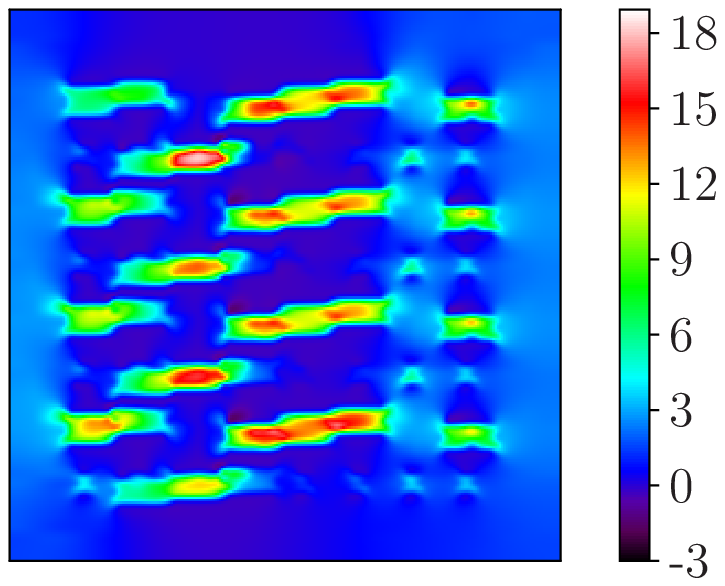}
    }
    \hspace*{-0.04\textwidth}
  \subfigure[Coarse $N_c = 647$]{\label{fig:fluxc6}
     \includegraphics[width = 0.33\textwidth, keepaspectratio = true]{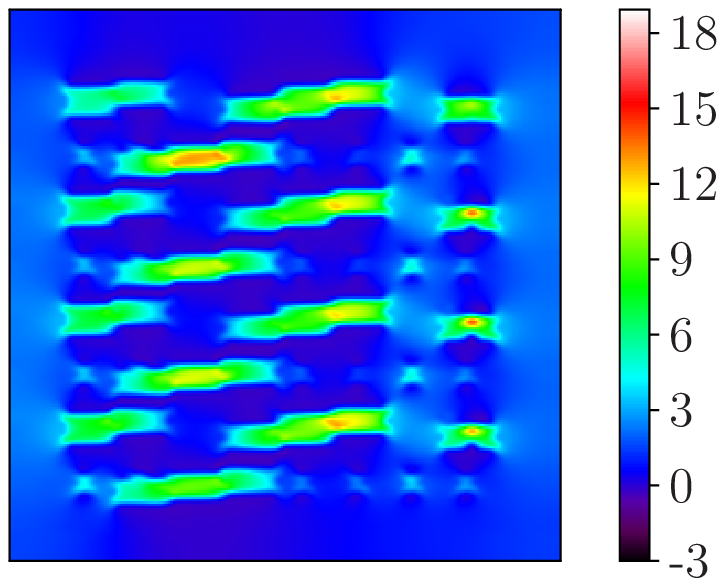}
    }
  
 \caption{ A comparison between the horizontal flux components of the fine system and the downscaled coarse system  }
 \label{fluxplot} 
\end{figure}

For more rigorous comparisons we also consider the relative error quantities given by 

\begin{equation} \label{error}
E_c =  \frac { \| p - p_c \| } { \| p \| } \times 100 \%,
\end{equation}
where $p$ denotes the reference fine-grid solution, $p_c$ is a specified coarse-grid (GMsFEM-FV) solution. For the norm quantities in \eqref{error}, we consider both the energy and weighted $L^2$ norms respectively given by 

$$\displaystyle \| p \|_{H^1_k (\Omega)} = \left( \int_\Omega k  | \nabla p |^2 \right)^{1/2} \quad \text{and} \quad \| p \|_{L^2_k (\Omega)} = \left( \int_\Omega k  p^2 \right)^{1/2}. $$ 
We note that the energy error is of particular importance for the analysis associated with GMsFEM (see, e.g., \cite{egh12}), and involves the gradient of the solution (which offers a flux-
like comparison). However, using either norm we are primarily interested in illustrating the effects of larger coarse spaces (i.e., more basis functions) in the GMsFEM-FV construction. In 
Table \ref{taberrors} we offer a number of errors corresponding to a variety of coarse space dimensions. The left most column tabulates the total size of the full coarse system 
associated with a specified level of enrichment and conservation constraints. The next two columns itemize the dimension of the coarse space (and corresponding number of basis 
functions per coarse node), and the number of FV constraints. Most importantly, the results indicate that an increase in the coarse-space dimension yields a predictable error decline 
associated with the solution and its gradient. As such, these results strongly suggest that incorporating numerous pressure solves into the context of the operator splitting technique 
(recall Section \ref{solalg}) will yield increasingly accurate two-phase saturation solutions. This is what we consider in the next subsection.  

\begin{table}
  \centering
  \begin{tabular}{ | c | c | c | c | c |}
    \hline
    Full System & Coarse Dimension  & Constraints & \multicolumn{2}{ | c |}{Relative Errors (\%)} \\
    \hline \hline
    $N_c$  & $\text{dim}(V_0)$  [$\#$ basis] & $M_c$ & $L^2_k (\Omega)$ & $H^1_k (\Omega)$ \\
    \hline
    242  & 121  [1]   &   121 & 8.4     &   $>$100   \\
    \hline
    323  & 202  [2]   &   121 & 7.4     &   39.8   \\
    \hline
    485  & 364  [4]  &   121 & 1.2     &   15.0    \\
    \hline
    647  & 526  [6]  &   121 & 0.9     &   11.1   \\
    \hline
    809  & 688  [8]  &   121 & 0.3     &   9.0     \\
    \hline
    971  & 850  [10]  &   121 & 0.3     &   8.1   \\
    \hline
  \end{tabular}
\caption{Relative pressure and flux errors for a variety of coarse space dimensions}
\label{taberrors}
\end{table}

\subsection{Two-phase saturation}
In this section we apply the constrained GMsFEM method to the full two-phase model as described in Section \ref{tpmodel}. More specifically, we apply the conservative GMsFEM 
discretization of Eq. \eqref{eq:problem1} for each pressure update. Then, a fine scale conservative flux field is obtained via the downscaling procedure in Section \ref{downscale} in 
order to march the saturation solution in time using the explicit scheme from Section \ref{solalg}. In doing so, we assess the effectiveness of the proposed approach in a context where 
it is repeatedly used to accurately capture a number of pressure equation updates. 

\begin{figure}[tb]
 \centering
   \includegraphics[width = 0.9\textwidth, keepaspectratio = true]{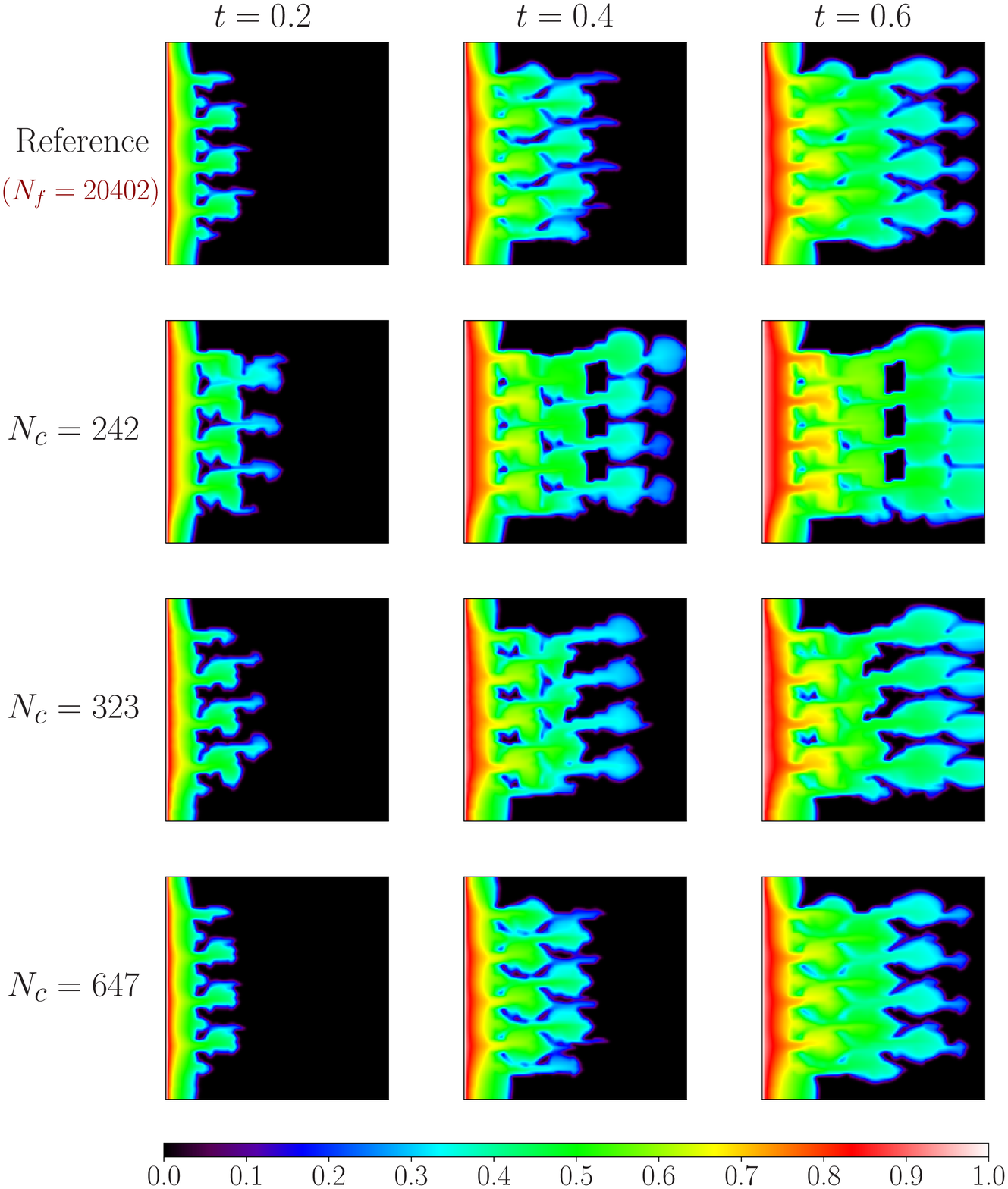}
 \caption{ Water saturation profiles advancing in time for a variety of reduced-order dimensions }
 \label{plotsat}
\end{figure}

To solve the two-phase model given in Eqs. \eqref{pres} and \eqref{sat} we use quadratic relative permeability curves given by $k_{rw} = S^2$ and $k_{ro} = (1 - S)^2$. We 
additionally use values of $\mu_w = 1$ and $\mu_o = 5$ for the water and oil phase viscosities. The same pressure boundary conditions from Section \ref{spresult} are used, and for 
the initial saturation condition, we set $S=1$ at the left edge and assume $S(x,0) = 0$ elsewhere. Finally, we recall that the high-contrast permeability coefficient illustrated in Fig. 
\ref{coeff} is used for $k (x)$. 
For a motivating application of the method, we offer a representative set of two-phase flow solutions in Fig. \ref{plotsat}. The plot shows saturation solutions advancing in time for a 
variety of coarse space dimensions, compared with fine-scale reference solutions. In particular, the first row shows three saturation snapshots for the case when $N_f = 20402$, and the 
second, third and fourth rows respectively shows saturation snapshots for the cases when $N_c = 242, 323, 647$. We note a significant improvement for the case when $N_c = 647$ as 
compared to the lower dimension $N_c = 121$. More specifically, the addition of more basis functions to the coarse pressure space yields flux and saturation values that accurately 
capture the fine-scale behavior of the system. And, with respect to the reduced dimension $N_c = 647$, we can see that the solutions are nearly indistinguishable from the reference 
solutions.

\begin{figure}[tb]
 \centering
   \includegraphics[width = 0.75\textwidth, keepaspectratio = true]{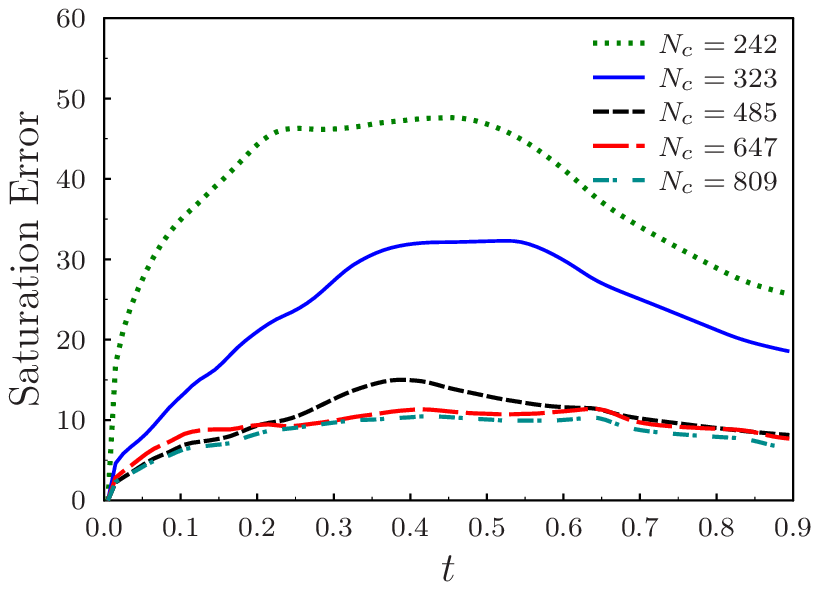}
 \caption{ $L^2$ error results for the two-phase flow problem }
 \label{twophase}
\end{figure}

For a final set of comparisons, we tabulate the standard $L^2$ relative error of the saturation profiles from Fig. \ref{plotsat}, as well as for a variety of other scenarios. The full set of 
results can be seen in Fig. \ref{twophase}. For these examples, we use a variety of coarse-space dimensions $N_c = 242, 323, 485, 647,$ and $809$ and run the simulations to a final 
time of $T = 0.9$. A time stepping value of $\Delta t = 10^{-4}$ is used, and we update the saturation for $100$ time steps in between each pressure solve. As was evident from the 
previous illustration, we can see from Fig. \ref{twophase} that the error may be significantly decreased by adding more basis functions in the GMsFEM construction. As a particular 
example, a maximum error value of roughly $50 \%$ ($N_c = 242$) may be decreased to a maximum error value of roughly $8 \%$ in the case when $N_c = 809$. These results 
serve to further illustrate the flexibility and accuracy associated with the conservative GMsFEM construction. More specifically, more basis functions may be used to obtain a higher level 
of accuracy, whereas less basis functions may be used when efficiency is a main consideration.

\section{Concluding remarks} \label{conclusion}
In this paper, we propose a method for the construction of locally conservative flux fields through a constrained variation of the Generalized Multiscale Finite Element Method (GMsFEM). 
The flux values are obtained through the use of a Ritz formulation in which we augment
the resulting linear system of the continuous Galerkin (CG) formulation in the higher-order GMsFEM
approximation space. We impose the finite volume-based restrictions through incorporating a scalar Lagrange multiplier for each mass conservation constraint on a specified scale, and 
as such, the proposed method may be viewed as a minimization problem in which an energy functional of the governing equations is minimized within a subspace of functions that 
satisfy the desired conservation properties. Due to the inherent construction of the coarse-grid solution space, and the way in which the mass conservation constraints are imposed, the 
combined methodology is shown to offer a robust and flexible framework for obtaining conservative flux fields to be used in two-phase flow modeling. 
To illustrate the performance of the method we consider model flow equations with heterogeneous permeability coefficients that have high-variation and discontinuities which 
significantly affect the flow patterns of the two-phase model. The increase in accuracy associated with the computation of the GMsFEM pressure solutions is inherited by the flux fields 
and saturation solutions, and is closely correlated to the size of the reduced-order systems. In particular, the addition of more basis functions to the enriched multiscale space produces 
solutions that more accurately capture the behavior of the fine scale model. A variety of single- and two-phase numerical examples are offered to validate the performance of the 
method.

\section*{Acknowledgments}
J. Galvis would like to thank R. Lazarov and P. Chatzipantelidis for interesting discussions on higher order  finite volume methods and for pointing out some of the references. M. Presho 
is supported by the U.S. Department of Energy Office of Science, Office of Advanced Scientific Computing Research, Applied Mathematics program under Award Number DE-SC0009286 
as part of the DiaMonD Multifaceted Mathematics Integrated Capability Center.

\section*{References}


\end{document}